\documentclass[11pt]{amsart}
\usepackage {amssymb, amsmath, amsgen, amstext, 
amsbsy, amsopn, amsfonts, amsthm}

\newcommand{\DynkinG}[1]{
\setlength{\unitlength}{#1 true in}

\begin{picture}(3.6, 1)

\multiput(1.2,.1)(1,0){2}{\circle*{.2}}
\put(1.2,.04){\line(1,0){1}}
\put(1.2,.1){\line(1,0){1}}
\put(1.2,.16){\line(1,0){1}}
\put(1.4,-.05){$>$}
\end{picture}}

\newcommand{\weightDynkinG}[3]{
\setlength{\unitlength}{#1 true in}

\begin{picture}(3.5, 1)

\put(.8,0){$#2$}
\put(1.8,0){$#3$}

\end{picture}}

\newcommand{\DynkinF}[5]{
\setlength{\unitlength}{#1 true in}

\begin{picture}(3.6, 1)

\put(.2,.1){\line(1,0){1}}
\multiput(.2,.1)(1,0){4}{\circle*{.2}}
\put(2.2,.1){\line(1,0){1}}
\put(1.2,.06){\line(1,0){1}}
\put(1.2,.14){\line(1,0){1}}
\put(1.43,-.05){$>$}
\put(.08,.27){$#2$}
\put(1.08,.27){$#3$}
\put(2.08,.27){$#4$}
\put(3.08,.27){$#5$}
\end{picture}}

\newcommand{\weightDynkinF}[5]{
\setlength{\unitlength}{#1 true in}

\begin{picture}(3.6, 1)

\put(0.08,0){$#2$}
\put(1.08,0){$#3$}
\put(2.08,0){$#4$}
\put(3.08,0){$#5$}
\end{picture}}

\newcommand{\DynkinEsix}[1]{
\setlength{\unitlength}{#1 true in}

\begin{picture}(4.6, 1)(0,-.55)
\multiput(.2,.1)(1,0){5}{\circle*{.2}}
\put(.2,.1){\line(1,0){4}}
\put(2.2,-.5){\circle*{.2}}
\put(2.2,.1){\line(0,-1){.5}}
\end{picture}}

\newcommand{\weightDynkinEsix}[7]{
\setlength{\unitlength}{#1 true in}

\begin{picture}(4.7, .9)
\put(0.08,0){$#2$}
\put(1.08,0){$#3$}
\put(2.08,0){$#4$}
\put(3.08,0){$#5$}
\put(4.08,0){$#6$}
\put(2.08,-.8){$#7$}

\end{picture}}

\newcommand{\DynkinEseven}[1]{
\setlength{\unitlength}{#1 true in}

\begin{picture}(5.7, 1)(0,-.55)
\multiput(.2,.1)(1,0){6}{\circle*{.2}}
\put(.2,.1){\line(1,0){5}}
\put(2.2,-.5){\circle*{.2}}
\put(2.2,.1){\line(0,-1){.5}}
\end{picture}}

\newcommand{\weightDynkinEseven}[8]{
\setlength{\unitlength}{#1 true in}

\begin{picture}(5.8, .9)
\put(0.08,0){$#2$}
\put(1.08,0){$#3$}
\put(2.08,0){$#4$}
\put(3.08,0){$#5$}
\put(4.08,0){$#6$}
\put(5.08,0){$#7$}
\put(2.08,-.8){$#8$}

\end{picture}}

\newcommand{\DynkinEeight}[1]{
\setlength{\unitlength}{#1 true in}

\begin{picture}(6.7, 1)(0,-.55)
\multiput(.2,.1)(1,0){7}{\circle*{.2}}
\put(.2,.1){\line(1,0){6}}
\put(2.2,-.5){\circle*{.2}}
\put(2.2,.1){\line(0,-1){.5}}
\end{picture}}

\newcommand{\weightDynkinEeight}[9]{
\setlength{\unitlength}{#1 true in}

\begin{picture}(6.8, .9)
\put(0.07,0){$#2$}
\put(1.07,0){$#3$}
\put(2.07,0){$#4$}
\put(3.07,0){$#5$}
\put(4.07,0){$#6$}
\put(5.07,0){$#7$}
\put(6.07,0){$#8$}
\put(2.07,-.8){$#9$}

\end{picture}}

\newcommand{\g}{{\mathfrak g}}

\newcommand{\levi}{{\mathfrak l}}

\newcommand{\flag}{{\mathcal B}}
\newcommand{\orbit}{{\mathcal O}}

\newcommand{\nilcone}{{\mathcal N}}
\newcommand{\dualnil}{{^L \mathcal N}}
\newcommand{\spnilcone}{{\mathcal N}^{sp}_{o}}
\newcommand{\dualg}{{ ^L \mathfrak g }}
\newcommand{\dualG}{{^L G}}

\newcommand{\ro}{{\Phi}}

\newcommand{\complex}{{\mathbf C}}

\newcommand{\nat}{\mathbb N}

\newcommand{\bval}{\tilde{b}}
\newcommand{\parti}{\mathcal P}

\newtheorem{proposition}{Proposition}

\newtheorem{lemma}[proposition]{Lemma}
\newtheorem{theorem}[proposition]{Theorem}

\theoremstyle{definition}

\newtheorem{remark}[proposition]{Remark}

\begin{document}

\setlength{\oddsidemargin}{.3in}
\setlength{\evensidemargin}{0in} 

\title { Lusztig's Canonical Quotient and Generalized Duality}
\author{ Eric Sommers}

\address{
Department of Mathematics \\
Harvard University \\
One Oxford St. \\
Cambridge, MA 02138 \\
U.S.A.}

\date{4/15/01}

\email{esommers@math.harvard.edu}

\begin{abstract}

We give a new characterization of Lusztig's canonical quotient, a finite group attached
to each special nilpotent orbit of a complex semisimple Lie algebra.  
This group plays an important role
in the classification of unipotent representations of finite groups of Lie type.

We also define a duality map.  To each pair of a nilpotent orbit and a conjugacy
class in its fundamental group, the map assigns a nilpotent orbit 
in the Langlands dual Lie algebra.
This
map is surjective and is related to a map introduced by Lusztig (and studied by 
Spaltenstein).  When the conjugacy
class is trivial, our duality map is just the one studied by Spaltenstein 
and by Barbasch and Vogan which has 
image consisting of the special nilpotent orbits.

\end{abstract}

\maketitle

\section{Introduction}

To each special nilpotent orbit in a simple Lie algebra over the complex numbers, Lusztig 
has assigned a finite group (called Lusztig's canonical quotient) which is naturally
 a quotient of the fundamental group of the orbit.  
This group is 
defined using the Springer correspondence and 
the generic degrees of Hecke algebra representations; it plays an
important role in Lusztig's work on parametrizing unipotent representations of a finite group
of Lie type \cite{lu3}. 

In this paper we give a description of the canonical quotient which uses 
the description of the component group of a nilpotent element 
given in \cite{so1} when $G$ is of adjoint type.
To do this we assign 
to each conjugacy class in the component group a numerical value (called
 the $\bval$-value).  The definition of the $\bval$-value 
makes use of the duality map of \cite{spa1}  
and involves the dimension of a Springer fiber; it 
is related to the usual $b$-value of a Springer representation
(the smallest degree in which the representation appears in the 
harmonic polynomials on the Cartan subalgebra). 
Our observation is that 
the kernel of the quotient map 
from the component group to the canonical quotient
consists of all conjugacy classes whose $\bval$-value is equal to 
the $\bval$-value of the trivial conjugacy class.

In the second part of the paper 
we define a surjective map from the set of pairs consisting of 
a nilpotent orbit 
and a conjugacy class in its fundamental group to the
set of nilpotent orbits in the Langlands dual Lie algebra.
This extends the duality of \cite{spa1} (which corresponds to the 
situation when the conjugacy class is
trivial).  In fact the map depends 
only on the image of the conjugacy class in
Lusztig's canonical quotient.  
Our map is defined by combining the duality of \cite{spa1}
and a map defined by Lusztig \cite{lu3}.  
Lusztig's map can be thought
of as a generalization of induction \cite{lusp}; 
in his book, Spaltenstein also studied
this generalization of induction \cite{spa1}.  

I would like to thank G. Lusztig, D. Vogan, P. Achar, P. Trapa,
and D. Barbasch for helpful conversations.  
I gratefully acknowledge the support of NSF grant DMS-0070674
and the hospitality of H. Kraft and the Mathematisches Institut
at the Universit\"at Basel in June, 1999.

\section{Basic set-up}

Throughout the paper, $G$ is a connected simple 
algebraic group over the complex numbers $\complex$.  
In $G$ we fix a maximal torus $T$ contained in a Borel subgroup $B$.
Let the character group of $T$ be $X^*(T)$. 
We use $\g$ for the Lie algebra of $G$. 
Let $\ro$ be the roots of $G$ and 
let $W = N_G(T)/ T$ be the Weyl group of $G$ (with respect to $T$).

The group $G$ acts on $\g$ via the adjoint action.  If $e \in \g$ is an element of
$\g$, we denote by $\orbit_e$ the orbit of $\g$ through $e$ under the action of $G$.
If $e$ is a nilpotent element, we call $\orbit_e$ a nilpotent orbit.  

We denote by $\nilcone_o$ the set of nilpotent
orbits in $\g$. This set is partially ordered by the relation
$\orbit_e \preceq \orbit_f$ whenever $\bar{\orbit}_e \subset \bar{\orbit}_f$.

For $\orbit \in \nilcone_o$ and $e \in \orbit$, let $A(e)$ be the
component group $Z_G(e)/ Z^{0}_G(e)$.  
If $e, e' \in \orbit$, then 
$A(e)$ may be identified with $A(e')$ and so we write
$A(\orbit)$ for this finite group.  
Furthermore, we may speak in a well-defined manner 
of the conjugacy classes of $A(\orbit)$.
When $G$ is simply-connected (and we pass to the analytic topology), 
$A(\orbit)$ is just the fundamental group of $\orbit$.
We denote by $\hat{A}(\orbit)$ the irreducible representations of 
$A(\orbit)$ and 
use the notion of irreducible local system on $\orbit$ interchangeably with
the notion of irreducible representation of $A(\orbit)$, 
or $A(e)$ for any $e \in \orbit$. 

We denote by $\nilcone_{o,c}$ the set of pairs $(\orbit, C)$ 
consisting of an orbit $\orbit \in \nilcone_o$ 
and a conjugacy class $C \subset A(\orbit)$. 
The set of special orbits in $\nilcone_o$ 
(see \cite{lu1}) will be denoted by 
$\spnilcone$.  
There is an order-reversing duality map $d:\nilcone_o \to \spnilcone$ 
studied by Spaltenstein such that $d^2$ is the identity on $\spnilcone$
\cite{spa1}. This map is already implicit in \cite{lu1}, so we refer
to $d$ henceforth as Lusztig-Spaltenstein duality.

\section{Notation on partitions}

In the classical groups it will be helpful  
to have a description of 
the elements of $\nilcone_o$ and $\nilcone_{o,c}$
and the map $d$
in terms of partitions.  We introduce that notation now
(roughly following the references \cite{cm}, \cite{ca}, and \cite{spa1}).

Let $\parti (m)$ denote the set of
partitions $\lambda = [\lambda_1 \geq \dots \geq \lambda_k]$ of $m$ 
(we assume that $\lambda_k \neq 0$).
For a part $\lambda_i$ of $\lambda$, we call $i$ the position of
$\lambda_i$ in $\lambda$.
For $\lambda \in \parti (m)$, let $r_i = \# \{j \ | \ \lambda_j = i \}$.
For $\epsilon \in \{ 0,1 \}$ let $\parti_{\epsilon} (m)$ be 
the set of partitions $\lambda$ 
of $m$ such that $r_i \equiv 0$ whenever $i \equiv \epsilon$ 
(all congruences are modulo $2$).  

It is well known that $\nilcone_{o}$ is in bijection with 
$\parti (n+1)$ when $\g$ is of type $A_n$; with $\parti_{1} (2n)$ when $\g$ is of type $C_n$;
with $\parti_{0} (2n+1)$ when $\g$ is of type $B_n$; and with  
$\parti_{0} (2n)$ when $\g$ is of type $D_n$, except that those partitions
with all even parts correspond to two orbits in $\nilcone_{o}$ (called
the very even orbits).
In what follows we will never have a need to address separately the
very even orbits,
so we will not 
bother to introduce notation to distinguish between
very even orbits.

We will also refer to $\parti_{1} (2n)$ as $\parti_{C} (2n)$;
to $\parti_{0} (2n+1)$ as  $\parti_{B} (2n+1)$;
and to $\parti_{0} (2n)$ as  $\parti_{D} (2n)$.

The set $\parti(m)$ 
is partially ordered by the usual partial ordering on partitions.
This induces a partial ordering on the sets  
$\parti_{C} (2n)$, $\parti_{B} (2n+1)$, and $\parti_{D} (2n)$
and these partial orderings coincide
with the partial ordering on nilpotent orbits given
by inclusion of closures.  We will refer to nilpotent orbits and partitions
interchangeable in the classical groups (with the
caveat mentioned earlier for the very even orbits in type $D$).

We can also represent elements of 
$\nilcone_{o,c}$ in terms of partitions
in the classical groups (see the last section of \cite{so1}).  
Given  $(\orbit, C) \in \nilcone_{o,c}$ let $e \in \orbit$
and let $s \in Z_G(e)$ be a semisimple element
whose image in $A(e) \cong A(\orbit)$ belongs to $C$.
Set $\levi = Z_{\g}(s)$.  
We may as well assume that $\levi$ contains Lie($T$); then 
the subalgebra $\levi$ of $\g$ 
(which we called a pseudo-Levi subalgebra in \cite{so1})
corresponds to a proper subset of the extended Dynkin diagram of $\g$.
It is always possible to choose 
$s$ so that $\levi$ has semisimple rank equal to $\g$,
which we will do.
Next we write $\levi = \levi_1 \oplus \levi_2$, where
$\levi_2$ is a semisimple subalgebra 
of the same type as $\g$
and $\levi_1$ is a simple Lie algebra (or possibly zero)
which contains
the root space corresponding to the lowest root of $\g$
(the extra node in the extended Dynkin diagram).
Then $\levi_1$ (if non-zero) 
is of type $D, C,$ or $D$ depending on whether $\g$ is 
of type $B, C,$ or $D$, respectively.
Write $e = e_1 + e_2$ where $e_1 \in \levi_1$ and $e_2 \in \levi_2$.
Finally, we may modify the choice of $s$ so that $e_1$ is a
distinguished nilpotent element in $\levi_1$.
Then we can attach to $(\orbit, C) \in \nilcone_{o,c}$
a pair of partitions $(\nu, \eta)$ where
$\nu$ is the partition of $e_1$ in $\levi_1$
and $\eta$ is the partition of $e_2$ in $\levi_2$.
We have
\begin{equation} \label{pairs}
\begin{aligned}
\nu \in \parti_D (2k), \ & \eta \in \parti_B (2n+1-2k) & \quad \text{ in type $B_n$} \\
\nu \in \parti_C (2k), \ & \eta \in \parti_C (2n-2k) &  \quad \text{ in type $C_n$} \\
\nu \in \parti_D (2k), \ & \eta \in \parti_D (2n-2k) &  \quad \text{ in type $D_n$} \\
\end{aligned}
\end{equation}
where $0 \leq k \leq n$ and 
$\nu$ is a distinguished partition,
meaning $r_i \leq 1$ for all $i$ 
(this necessarily forces 
$r_i = 0$ whenever $i \equiv \epsilon$).
The partition $\lambda$ of $\orbit$ 
is just the partition consisting of all the parts in $\nu$ and $\eta$.
We denote this partition by $\nu \cup \eta$ and write $\lambda = \nu \cup \eta$.
In types $C$ and $D$ we make the further assumption 
that the largest $i \not\equiv \epsilon$ with $r_i \equiv 1$
is not a part of $\nu$.  These assumptions guarantee that 
when $G$ is of adjoint type there is a bijection between 
$\nilcone_{o,c}$ and the pairs $(\nu, \eta)$ specified above
(if $\orbit$ is very even in type $D$, the component group
$A(\orbit)$ is trivial when $G$ is adjoint; 
so there are no additional complications beyond the one mentioned earlier). 
We note that $\nu$ will be the empty partition (that is, $k=0$ in 
equation ~(\ref{pairs})) 
if and only if $C$
is the trivial conjugacy class.

For example, let $\orbit$ be the orbit with partition $[5,3,1]$ in $B_4$.
Then we have 
$$(\emptyset, [5,3,1]), ([3,1], [5]), ([5,1],[3]),([5,3], [1])$$
are the four elements of $\nilcone_{o,c}$ corresponding to the four
conjugacy classes of $A(\orbit)$ when $G$ is of adjoint type,
with the first pair corresponding to the trivial conjugacy class.

The dual partition of $\lambda \in \parti (m)$, denoted 
$\lambda^*$, is defined
by $\lambda^*_i = \# \{ j \ | \ \lambda_j \geq i \}$.
Let $X=B$, $C$, or $D$.
We define the $X$-collapse of a partition $\lambda \in \parti(m)$
where $m$ is even when $X = C$ or $D$ and $m$ is odd when $X=B$.  
The $X$-collapse of $\lambda \in \parti(m)$ is the
partition $\mu \in \parti_{X}(m)$ such that 
$\mu \preceq \lambda$ and such that $\mu' \preceq \mu$
whenever $\mu' \in \parti_{X}(m)$ and
$\mu' \preceq \lambda$.
The $X$-collapse of $\lambda$ is denoted $\lambda_X$.
It is well-defined and unique. 

The duality map $d:\nilcone_o \to \spnilcone$ 
can now be expressed as follows in the classical groups.
In type $A$, we have $d( \lambda) = \lambda^*$
and in type $X$ (for $X= B, C, D$) 
we have  $d( \lambda) = (\lambda^*)_X$.
The set of special nilpotent orbits correspond to the partitions in the
image of $d$ (in type $D$ all very even orbits 
are special).  
Hence in type $A$ all nilpotent orbits are special.
In the other types it is known that
$\lambda$ is special if and only if
\begin{equation} 
\begin{aligned}
\lambda^* &\in \parti_B(2n+1) \text{ in type $B_n$} \\
\lambda^* &\in \parti_C(2n) \text{ in type $C_n$} \\
\lambda^* &\in \parti_C(2n) \text{ in type $D_n$}.  
\end{aligned}
\end{equation}

\section{$\bval$-value}

In this section we assign to each pair $(\orbit, C) \in \nilcone_{o,c}$ 
a natural number which we will call the $\bval$-value of $(\orbit, C)$ and
which will be denoted by $\bval_{(\orbit, C)}$.  

Given  $(\orbit, C) \in \nilcone_{o,c}$ let $e \in \orbit$ 
and $s \in Z_G(e)$ be a semisimple element
whose image in $A(e) \cong A(\orbit)$ belongs to $C$.
The group $L= Z_G(s)$ has as its Lie algebra $\levi = Z_{\g}(s)$ and clearly
$e \in \levi$.

Let $\orbit^{\levi}_e$ denote the orbit in $\levi$
through $e$ under the action of $L$.
Applying Lusztig-Spaltenstein duality $d_{\levi}$ with respect to $\levi$, 
we obtain the orbit
$\orbit^{\levi}_f= d_{\levi}(\orbit^{\levi}_e) \subset \levi$. 
We define 
$\bval_{(\orbit, C)} = 
\frac{1}{2}( \dim (\levi) - \dim(\orbit^{\levi}_f) - \dim (T))$.
Equivalently, if 
$\flag^{\levi}_f$ is the variety of Borel subalgebras of $\levi$ 
which contain $f$, then $\bval_{(\orbit, C)} = \dim \flag^{\levi}_f$.
We have the following proposition whose proof we give after the proof
of proposition \ref{indep}.
\begin{proposition} \label{keyprop}

The number $\bval_{(\orbit, C)}$ is well-defined.  
In other words, it is independent of the choices made for $s$ and $e$.

\end{proposition}

By results in \cite{so1} it is always possible to choose $s$ so that
$e$ is distinguished in $\levi = Z_{\g}(s)$.
By the previous proposition,  
in order to compute $\bval_{(\orbit, C)}$ it is enough
to determine a pair $(\levi, e)$
attached to $(\orbit, C)$ where $e$ is distinguished in $\levi$
(this is done in \cite{so1}) and to compute the 
$\bval$-value of $e$ with respect to $\levi$ (with the trivial conjugacy class
$1$ of $A(e)$).
Hence it suffices to compute 
$\bval_{(\orbit, 1)}$ for each distinguished orbit $\orbit$ in 
each simple Lie algebra.

We now record 
$\bval_{(\orbit, 1)}$ for the distinguished orbits
in the exceptional groups.

\bigskip
\begin{center}

\begin{tabular}{|c|c|c|} \hline
\multicolumn{3}{|c|}{$G_{2}$} \\ \hline
\multicolumn{1}{|r|}{$\DynkinG{.2}$} 
& \multicolumn{1}{|c|}{Bala-Carter }
& \multicolumn{1}{|c|}{$\bval$} \\ \hline
\weightDynkinG{.18}{2}{2} &  $G_{2}$  & 6\\ \hline
\weightDynkinG{.18}{2}{0} &  $G_{2}(a_{1})$  & 1\\ \hline
\end{tabular}

\medskip

\begin{tabular}{|c|c|c|} \hline
\multicolumn{3}{|c|}{$F_{4}$} \\ \hline
\multicolumn{1}{|r|}{$\DynkinF{.2}{}{}{}{}$} 
& \multicolumn{1}{|c|}{Bala-Carter }
& \multicolumn{1}{|c|}{$\bval$} \\ \hline
\weightDynkinF{.2}{2}{2}{2}{2} & $F_4$ & 24 \\ \hline
\weightDynkinF{.2}{2}{2}{0}{2} & $F_4(a_1)$ &  13 \\  \hline
\weightDynkinF{.2}{0}{2}{0}{2} & $F_4(a_2)$ &  10 \\ \hline
\weightDynkinF{.2}{0}{2}{0}{0} & $F_4(a_3)$ &  4 \\ \hline

\end{tabular}

\medskip

\begin{tabular}{|c|c|c|} \hline
\multicolumn{3}{|c|}{$E_{6}$} \\ \hline
\multicolumn{1}{|c|}{\DynkinEsix{.18}} 
& \multicolumn{1}{|c|}{Bala-Carter }
& \multicolumn{1}{|c|}{$\bval$} \\ \hline
\weightDynkinEsix{.18}{2}{2}{2}{2}{2}{2} & $E_6$ & 36 \\  
 & & \\ \hline
\weightDynkinEsix{.18}{2}{2}{0}{2}{2}{2} & $E_6(a_1)$ & 25 \\  & & \\  \hline
\weightDynkinEsix{.18}{2}{0}{2}{0}{2}{0} & $E_6(a_3)$ & 15 \\  & & \\  \hline

\end{tabular}

\medskip

\begin{tabular}{|c|c|c|} \hline
\multicolumn{3}{|c|}{$E_{7}$} \\ \hline
\multicolumn{1}{|c|}{\DynkinEseven{.18}} 
& \multicolumn{1}{|c|}{Bala-Carter }
& \multicolumn{1}{|c|}{$\bval$} \\ \hline
\weightDynkinEseven{.18}{2}{2}{2}{2}{2}{2}{2} & $E_7$ & 63 \\  & & \\ \hline
\weightDynkinEseven{.18}{2}{2}{0}{2}{2}{2}{2} & $E_7(a_1)$ & 46 \\  & & \\ \hline
\weightDynkinEseven{.18}{2}{2}{0}{2}{0}{2}{2} & $E_7(a_2)$ & 37 \\  & & \\ \hline
\weightDynkinEseven{.18}{2}{0}{2}{0}{2}{2}{0} & $E_7(a_3)$ & 30 \\ & & \\  \hline
\weightDynkinEseven{.18}{2}{0}{2}{0}{0}{2}{0} & $E_7(a_4)$ & 22 \\  & & \\ \hline
\weightDynkinEseven{.18}{0}{0}{2}{0}{0}{2}{0} & $E_7(a_5)$ & 16\\  & & \\ \hline

\end{tabular}

\bigskip

\begin{tabular}{|c|c|c|} \hline
\multicolumn{3}{|c|}{$E_{8}$} \\ \hline
\multicolumn{1}{|c|}{\DynkinEeight{.18}}
& \multicolumn{1}{|c|}{Bala-Carter }
& \multicolumn{1}{|c|}{$\bval$} \\ \hline
\weightDynkinEeight{.18}{2}{2}{2}{2}{2}{2}{2}{2} & $E_8$ & 120 \\  & & \\ \hline
\weightDynkinEeight{.18}{2}{2}{0}{2}{2}{2}{2}{2} & $E_8(a_1)$ & 91 \\ & & \\  \hline
\weightDynkinEeight{.18}{2}{2}{0}{2}{0}{2}{2}{2} & $E_8(a_2)$ & 74 \\  & & \\ \hline
\weightDynkinEeight{.18}{2}{0}{2}{0}{2}{2}{2}{0} & $E_8(a_3)$ & 63 \\ & & \\  \hline
\weightDynkinEeight{.18}{2}{0}{2}{0}{2}{0}{2}{0} & $E_8(a_4)$ & 52 \\  & & \\ \hline

\weightDynkinEeight{.18}{2}{0}{2}{0}{0}{2}{2}{0} & $E_8(b_4)$ & 47 \\ & & \\  \hline
\weightDynkinEeight{.18}{2}{0}{2}{0}{0}{2}{0}{0} & $E_8(a_5)$ & 42 \\  & & \\ \hline
\weightDynkinEeight{.18}{0}{0}{2}{0}{0}{2}{2}{0} & $E_8(b_5)$ & 37 \\ & & \\  \hline
\weightDynkinEeight{.18}{0}{0}{2}{0}{0}{2}{0}{0} & $E_8(a_6)$ & 32 \\  & & \\ \hline
\weightDynkinEeight{.18}{0}{0}{2}{0}{0}{0}{2}{0} & $E_8(b_6)$ & 28 \\ & & \\  \hline
\weightDynkinEeight{.18}{0}{0}{0}{2}{0}{0}{0}{0} & $E_8(a_7)$ & 16 \\ & & \\  \hline
\end{tabular}

\end{center}

\medskip

We now record $\bval_{(\orbit, 1)}$ for the distinguished orbits
in the classical groups.  In fact, it is no harder to record 
the formula for $\bval_{(\orbit, 1)}$ for any orbit $\orbit$, whether
distinguished or not.
Let $\lambda = [\lambda_1 \geq \lambda_2 \geq \dots \geq \lambda_k]$ 
be the partition for $\orbit$
in the appropriate classical group.  
Then $\bval_{(\orbit , 1)}$ equals
\begin{equation} 
\begin{aligned}
\frac{1}{2} \sum^{k}_{i=1} \lambda^2_i - 
\frac{1}{2} \sum^{k}_{i=1} \lambda_i  \text{ in type $A$} \\
\frac{1}{4} \sum^{k}_{i=1} \lambda^2_i - 
\frac{1}{2} \sum_{i \text{ \ odd}} \lambda_i + \frac{1}{4}
\text{ in type $B$} \\
\frac{1}{4} \sum^{k}_{i=1} \lambda^2_i - 
\frac{1}{2} \sum_{\text{$i$ \ even}} \lambda_i 
\text{ in types $C$ and $D$.}
\end{aligned}
\end{equation}
We omit the calculations (which are easy in all types
except type $D$) since they will follow quickly from later work. 

\begin{remark}
If $\orbit$ is the Richardson orbit for the 
parabolic subgroup $P$ 
whose Levi subgroup has root system  
isomorphic to $A_1 \times \dots \times A_1$ ($s$-times),
we have observed that 
$\bval_{(\orbit, 1)} = N -  s \ m_{n+1-s}$.  Here, $n$ is the rank of $G$,
$N$ is the number of positive roots of $G$, 
and $m_i$ is the $i$-th exponent of $G$ when the exponents are listed
in increasing order. 
This offers the hope that there is a way to determine the $\bval$-value
of any distinguished orbit from standard 
data arising purely from
the root system. 
\end{remark}


The next proposition is analogous to a statement (due to Lusztig in \cite{lu3})
about $a$-values of Springer representations
which we will recall in the next section. 

\begin{proposition} \label{inequality}

Let $1$ be the trivial conjugacy and let $C$ be any conjugacy class in $A(\orbit)$.
Then $\bval_{(\orbit, 1)} \leq \bval_{(\orbit, C)}$.  

\end{proposition}

\begin{proof}

In the exceptional groups, we have verified this directly.
In the classical groups, we must show that for any $\lambda \in \nilcone_o$
and $(\nu, \eta) \in \nilcone_{o,c}$ with $\lambda = \nu \cup \eta$
that we have 
$\bval_{\lambda} \leq \bval_{(\nu, \eta)} = 
\bval_{\nu} + \bval_{\eta}$ where the $\bval$-value is 
computed with respect to the appropriate subalgebras for $\lambda$, $\nu$,
and $\eta$.

In type $B$, $\bval_{\lambda} =\frac{1}{4} \sum^{k}_{i=1} \lambda^2_i - 
\frac{1}{2} \sum_{i \text{ \ odd}} \lambda_i + \frac{1}{4}$.
For any $(\nu, \eta)$ with $\lambda = \nu \cup \eta$, we must use 
both the $\bval$-value formula in type $B$ for $\eta$ and 
in type $D$ for $\nu$.  
It is clear that 
$\bval_{\nu} + \bval_{\eta} = \frac{1}{4} \sum^{k}_{i=1} \lambda^2_i
+ \frac{1}{4}$ minus various $\lambda_i$.  It turns out that 
that for each $m$ either $\lambda_{2m-1}$ or $\lambda_{2m}$ will be subtracted
in the formula, but not both, as we now show.
There are four possibilities for each $m$.  
If both $\lambda_{2m-1}$ and $\lambda_{2m}$ 
belong to $\eta$ or both belong to $\nu$,
then $\bval_{(\nu, \eta)}$ will have a term
$-\lambda_{2m-1}$ or $-\lambda_{2m}$, but not both, since 
the formulas for $B$ and $D$ select every other $\lambda_i$
to subtract.  If on the other hand,
$\lambda_{2m-1}$ belongs to $\eta$ and $\lambda_{2m}$ belongs to $\nu$, 
then again $\bval_{(\nu, \eta)}$ will have a term
$-\lambda_{2m-1}$ or $-\lambda_{2m}$, but not both.  
This is because the parity of the position of  
$\lambda_{2m-1}$ in $\eta$
will be the same as the parity of the position of 
$\lambda_{2m}$ in $\nu$
(as the position of $\lambda_{2m}$ in $\lambda$ is an even number, namely $2m$).  
But the formula for $B$ and $D$ choose to subtract
parts whose positions have opposite parity. Hence, only one can be selected
in the formula for $\bval_{\nu} + \bval_{\eta}$.  The result is the same
if $\lambda_{2m-1}$ belongs to $\nu$ and $\lambda_{2m}$ 
belongs to $\eta$.
Thus in either of the four cases, the effect is to subtract
a number which is less than or equal to $\lambda_{2m-1}$ since 
$\lambda_{2m-1} \geq \lambda_{2m}$, hence the inequality
$\bval_{\lambda} \leq \bval_{\nu} + \bval_{\eta}$.

In types $C$ and $D$, a similar argument holds, except we look
at the consecutive parts $\lambda_{2m}$ and $\lambda_{2m+1}$
of $\lambda$ (the part $\lambda_1$ will never be subtracted
in the formulas for the $\bval$-value).
One or the other, but not both, of these parts will be subtracted
in the formula for $\bval_{\nu} + \bval_{\eta}$
demonstrating the inequality.
\end{proof}

We recall some notation from \cite{so1} in the classical groups
\begin{gather*}
S_{odd} = \{ i \in \nat | \ i \not\equiv \varepsilon, \, r_i \equiv 1  \} \\
S_{even} = \{ i \in \nat | \ i \not\equiv \varepsilon, \, r_i \equiv 0 
\text{ and } r_i \neq 0 \}.
\end{gather*}
List the elements of $S_{odd}$ in decreasing order 
$j_l \geq j_{l-1} \geq \dots \geq j_2 \geq j_1$.
Assume that $l$ is even in type $C$ by 
setting $j_1=0$ if necessary
($l$ is automatically odd in type $B$ and automatically even in type $D$).

An element $(\nu, \eta) \in \nilcone_{o,c}$ determines
two sets $T_1$ and $T_2$ 
such that $T_{1} \subset S_{odd}$ and $T_{2} \subset S_{even}$
coming from the parts of $\nu$. 
Namely, the parts of $\nu$ (which each
occur with multiplicity one since $\nu$ is distinguished)
consist precisely of the elements in $T_1 \cup T_2$.  
For $j \in S_{odd}$, let 
$\delta_j = 1$ if $j \in T_1$ and $0$ if $j \notin T_1$.
Note we are assuming that $\delta_l = 0$ in types $C$ and $D$.
We define subsets $T^{(m)}_2$ of $T_2$ as follows:
let $T^{(m)}_2$ consist of those $i \in T_2$ 
such that $j_{m+1} \geq i \geq j_m$
and define $t_m$ to be the cardinality of $T^{(m)}_2$.

The next proposition follows easily from the previous
proposition and the formulas for the $\bval$-value.
\begin{proposition} \label{canonical}
In the classical groups
(not of type $A$),
the equality $\bval_{(\orbit, 1)} = \bval_{(\orbit, C)}$
holds if and only if
$t_m=0$ when $m$ is even
and $\delta_{m+1} + t_{m} + \delta_{m}$ is even when $m$ is odd.
\end{proposition}

\section{Lusztig's canonical quotient}

Our first main result is a new description of Lusztig's canonical quotient
$\bar{A}(\orbit)$, which is a quotient of $A(\orbit)$. 
These finite groups play an important role in Lusztig's classification of
unipotent representations of finite groups of Lie type.  Namely,
the set of unipotent representations for $G$ (split, with connected center, over a 
finite field) is parametrized by the following data:
a special nilpotent orbit $\orbit$, an element $x \in \bar{A}(\orbit)$,
and an irreducible representation of the centralizer of $x$ in 
$\bar{A}(\orbit)$ 
(all up to the appropriate conjugation).   
We hope to obtain a better 
understanding of this parametrization by having such an explicit description of the
canonical quotient and its conjugacy classes. 

Let us recall Lusztig's definition of $\bar{A}(\orbit)$.
Although Lusztig assumed $\orbit$ is special, 
his definition remains valid even if $\orbit$ is not special, 
so in what follows we do not assume $\orbit$ is special.

Recall that for each nilpotent element $e$ and local system $\pi$ 
on $\orbit_e$,
Springer has defined a representation $E_{e,\pi}$ which (if non-zero) is an irreducible
representation of $W$.  
Recall also that each irreducible representation $E$ of $W$
comes with two important
numerical invariants.   One comes from the generic degree of $E$ 
(the $a$-value) and one comes from the fake degree of $E$ (the $b$-value).  
We refer 
to \cite{lu3} for the definitions. 
Note that our notation is consistent with \cite{lu3}, but is not
consistent with \cite{lu1} or \cite{ca}.
In those sources, our $b$-value
is their $a$-value and our $a$-value is their $\tilde{a}$-value.  
The original definition of the canonical quotient of $A(\orbit)$ is as follows.  
Given $e \in \orbit$, consider the set 
$$\mathcal S = \{ \pi \in {\hat A(e)} \ | \  E_{e, \pi} \neq 0 \text{ and } a_{E_{e,\pi}} = a_{E_{e, 1}} \} $$ 
where ${\hat A(e)}$ is the set of irreducible representations of $A(e)$ 
and $E_{e,1}$ denotes the Springer representation associated to the 
trivial representation of $A(e)$.
Let $H$ be the intersection of
all the kernels of the representations in $\mathcal S$. 
Then $\bar{A}(\orbit)$ is defined to be the quotient  
$A(\orbit)/ H$ \cite{lu3}.
If a local system $\pi$ is not equivariant for $G$
when $G$ is of adjoint type, 
then $E_{e, \pi}$ will be zero.  Hence, the 
canonical quotient is the same for groups
in the same isogeny class of $G$
and so we assume $G$ is of adjoint type in what follows.   

\begin{remark}
Lusztig has observed that for 
$\pi \in {\hat A(e)}$ and $E_{e, \pi} \neq 0$
that $a_{E_{e,\pi}} \leq  a_{E_{e, 1}}$.  
Compare this with proposition \ref{inequality}.
\end{remark}

In light of Lusztig's original definition and proposition \ref{inequality}
we consider all conjugacy classes $C$ in $A(\orbit)$ with the property that
$\bval_{(\orbit, C)} = \bval_{(\orbit, 1)}$.  Let $H'$ be
the union of all such conjugacy classes.

\begin{theorem}
The set $H'$ coincides with $H$, so that
$\bar{A}(\orbit) = A(\orbit)/ H'$. 
Moreover, the  
$\bval$-value is constant on the cosets of $H'=H$ in $A(\orbit)$
(which are always a union of conjugacy classes). 
\end{theorem}

\begin{proof}
We are assuming that $G$ is of adjoint type (the result
is still valid for any $G$).
In the exceptional groups we verified the results directly
using the tables for $a$-values and Springer representations in \cite{ca}, 
the tables of conjugacy classes
in $A(\orbit)$ in \cite{so1}, and 
knowledge of the $\bval$-values for distinguished orbits given above.

In the classical groups we need to do some work in order to understand which local
systems appear in the set $\mathcal S$.

We illustrate the situation in type $B_n$.
Let $\orbit$ have partition $\lambda = [\lambda_1 \geq \dots \geq 
\lambda_k]$.  
As before, $r_i$ is the number of parts of $\lambda$ of size $i$.  
So $r_i$ is even whenever $i$ is even as we are in type $B$.
We associate to $\orbit$ a symbol as in chapter 13.3 of \cite{ca}.
It consists of $k$ elements.  
Let $d_i = \sum_{j<i} r_j$.
Each odd $i$ contributes to the symbol the interval of length $r_i$
$$(d_i + \frac{i-1}{2}, d_i + \frac{i-1}{2} + 1, \dots,  d_i + \frac{i-1}{2} + r_i - 1)$$
and each even $i$ contributes to the symbol the $r_i/2$ numbers 
$$d_i  + i/2, d_i + i/2 + 2, \dots,  d_i + i/2 + r_i - 2$$ 
each repeated twice in weakly increasing order.

Consider the elementary $2$-group with basis consisting of elements $x_i$, 
one for each 
odd $i$ with $r_i \neq 0$.  
Then $A(\orbit)$ is the subgroup of this group consisting of 
elements expressible as a sum of an even number of basis elements
(this is because we are working in the special
orthogonal group and not in the full orthogonal group).
Representations $\pi$ of $\hat{A}(\orbit)$ are thus specified 
by their values (of $\pm 1$)
on the $x_i$'s.
The representations $\pi$ for which $E_{e, \pi}$ is non-zero are those 
with the property that 
$$\# \{ x_i \ | \ r_i \equiv 1, d_i \equiv 0 \text{ and } \pi(x_i) = -1 \} =
\# \{ x_i \ | \ r_i \equiv 1, d_i \equiv 1 \text{ and } \pi(x_i) = -1 \}.$$
If $r_i$ is even, $\pi$ may have the value $1$ or $-1$ on $x_i$.

To determine the set $\mathcal S$ it is necessary 
to compute the $a$-value of each non-zero representation $E_{e, \pi}$
(see chapter 11.4 of \cite{ca}).  
We have to convert between two different notations for symbols.
This is a bit of a pain, but the work is greatly simplified
since we are only interested 
in when the $a$-value of $E_{e, \pi}$ equals the $a$-value of $E_{e, 1}$.
We find that the set $\mathcal S$ consists
of the following $\pi$.
Above we listed those odd $i$ with $r_i$ odd (which was denoted
$S_{odd}$ above) in decreasing order as 
$j_l \geq j_{l-1} \geq \dots \geq j_1$ (note that $l$ must be odd).  
If $\pi \in \mathcal S$, we must have $\pi(x_{j_l}) = 1$ and 
$\pi(x_{j_{2m}}) = \pi(x_{j_{2m-1}})$.
Let $i$ be odd and $r_i$ even.  
If $j_{2m+1} \geq i \geq j_{2m}$ for some $m$,
then $\pi$ may take values $\pm 1$ on $x_i$.
If, on the other hand, $j_{2m} \geq i \geq j_{2m-1}$, then 
$\pi(x_i) = \pi (x_{j_{2m-1}})$.

Since $A(\orbit)$ is abelian, each element forms its own conjugacy
class.  We need to relate our two descriptions of conjugacy
classes in $A(\orbit)$.
Given $x \in A(\orbit)$ write 
$x = x_{i_1} + x_{i_2} + \dots + x_{i_{2m}}$
where $i_1 \geq i_2 \geq \dots \geq i_{2m}$
(the usual classical description).
Then set $\nu = [i_1 \geq \dots \geq i_{2m}]$
and define $\eta$ by $\lambda = \nu \cup \eta$.
Then $(\nu, \eta) \in \nilcone_{o,c}$ exactly corresponds
to the conjugacy class of  $x \in A(\orbit)$ (see \cite{so2}).
Using our previous notation of $T^{(m)}_2$, $t_m$, and $\delta_j$,
we see that $H$ consists of those conjugacy classes
where $t_m=0$ when $m$ is even 
and $\delta_{m+1} + t_{m} + \delta_{m}$ is even when $m$ is odd.
By proposition \ref{canonical} 
this is exactly the condition that the conjugacy class
$(\nu, \eta)$ belongs to $H'$, showing that $H' = H$.
Finally, the statement that $\bval$ is constant on the cosets on 
$H'=H$ is an easy computation similar to the one done to 
find which classes belonged to $H'$.

A similar proof holds in type $C$ and $D$ which we omit.
\end{proof}

For completeness we record $\bar{A}(\orbit)$  
in the classical groups.  In the exceptional groups, $\bar{A}(\orbit)$
is listed in the last section.

In type $A$,  $\bar{A}(\orbit)$ is trivial.
In the other classical types, 
$\bar{A}(\orbit)$ is an elementary $2$-group.  Assuming $G$ is 
the special orthogonal group in types $B$ and $D$ and
$G$ is the symplectic group in type $C$,
we describe a subgroup $K$ of $A(\orbit)$ which 
maps bijectively onto $\bar{A}(\orbit) = A(\orbit)/H$.

In types $B$ and $D$, consider the subgroup $K$
of $A(\orbit)$ consisting of all elements expressible as 
a sum of an even number of $x_i$ where
$i$ is equal to some $j_m$ for $m$ odd 
or $i \in T^{(m)}_2$ for $m$ even.
In type $C$, consider the subgroup $K$
of $A(\orbit)$ consisting of all elements expressible as 
a sum of $x_i$
where $i$ is equal to some $j_m$ for $m$ odd 
or $i \in T^{(m)}_2$ for $m$ even.
Then $K \cong \bar{A}(\orbit)$.
In other words, the $i$ in question correspond, 
in type $B$, to corners of the Young diagram
of $\lambda$ which have odd length and odd height ($i$ is odd and $l$ is odd);
in type $D$, to corners of the Young diagram
of $\lambda$ which have odd length and even height ($i$ is odd and $l$ is even);
in type $C$, to corners of the Young diagram
of $\lambda$ which have even length and even height
($i$ is even and $l$ is even).

\section{Duality}

Let $\dualG$ be the Langlands dual group of $G$ with Lie algebra $\dualg$
(in other words, the root data of $\dualG$ and $G$ are dual).
Denote by $^L \nilcone_o$ the set of nilpotent orbits in $\dualg$
and $^L \nilcone^{sp}_o$ the set of special nilpotent orbits in $\dualg$.

There is a natural order-preserving (and dimension-preserving)
bijection between $\nilcone^{sp}_o$ and ${^L \nilcone^{sp}_o}$.
Therefore we can also view Lusztig-Spaltenstein duality 
as a map $d:\nilcone_o \to {^L \nilcone^{sp}_o}$.
It is in this context that Barbasch and Vogan have given a 
representation-theoretic description of Lusztig-Spaltenstein duality using primitive ideals \cite{bv}.
To distinguish between the duality
which stays within $G$ and the version which passes to the dual group,
we use the notation $d(\orbit)$ to mean the dual orbit within $\g$
and we use the notation 
$d_{(\orbit, 1)}$ to mean the dual orbit within $\dualg$. 
This is consistent with later (and earlier) notation.

We will now define a map $d: \nilcone_{o,c} \to \dualnil_{o}$ 
with the property that it
extends Lusztig-Spaltenstein duality.
It will turn out that $d$ is surjective and depends only 
on the image of $C$ in the canonical quotient ${\bar A}(\orbit)$.
Moreover, if $e' \in d_{(\orbit , C)}$, then we will have 
$\mbox{dim} ({^L \flag}_{e'}) = \bval_{(\orbit, C)}$
where $^L \flag$ denotes the flag variety of $\dualG$.

The definition of $d$ is as follows.
Given $(\orbit , C) \in \nilcone_{o,c}$, pick $e \in \orbit$
and a semisimple element $s \in Z_{G}(e)$
whose image in $A(e) \cong A(\orbit)$ belongs to $C$.
The group $L= Z_G(s)$ has as its Lie algebra $\levi = Z_{\g}(s)$ and clearly
$e \in \levi$.
Let $\orbit^{\levi}_e$ denote the orbit in $\levi$
through $e$ under the action of $L$.
Applying Lusztig-Spaltenstein duality $d_{\levi}$ with respect to $\levi$
(we assume here that $d_{\levi}$ stays within 
$\levi$ and we don't pass to the dual group), 
we obtain the orbit
$\orbit^{\levi}_f= d_{\levi}(\orbit^{\levi}_e) \subset \levi$. 
At this point we simply apply Lusztig's map from Chapter 13 of \cite{lu3}.

But first we recall the following result due to Joseph and 
proved uniformly in \cite{bm}.  
Every Weyl group representation of the form $E_{e,1}$
possesses property (B) of \cite{lusp} with respect to $G$.  That is,
the $b$-value of $E_{e,1}$ coincides with 
$\dim(\flag_e)$
and moreover the multiplicity of $E_{e,1}$ in the harmonic
polynomials of degree $\dim(\flag_e)$ on a Cartan subalgebra of $\g$
is exactly one.  

The Springer correspondence for $L$ 
produces an irreducible representation $E_{f, 1}$ (always non-zero)
of $W(s)$, the Weyl group
of  $\levi = Z_{\g}(s)$.  
The second part of the statement of property (B)
means that we can apply the operation of 
truncated induction to $E_{f,1}$ (see \cite{lu1}). 
So let $E = j_{W(s)}^W (E_{f, 1})$ be the representation of $W$ obtained by 
truncated induction of $E_{f, 1}$ from $W(s)$ to $W$.  
Then $E$ has the property that 
$b_E = b_{E_{f,1}}$ from the definition of truncated induction.  
Now consider $E$ as a representation of $^L W$, the Weyl group
of $\dualG$ (since $^L W$ is isomorphic to $W$ via the involution which
interchanges long and short roots). 
Lusztig has observed that 
$E$ is always of the form $E_{e', 1}$ for some nilpotent element $e'$ in 
$\dualg$ (this will be re-verified explicitly below).
We then define $d_{(\orbit , C)}$ to be the orbit in $\dualg$ through $e'$.  

It is clear from the definition (assuming it is well-defined)
and the first part of the statement of 
property (B) (applied twice, once in $L$ and
once in $\dualG$) 
that $\mbox{dim} ({^L \flag}_{e'}) = \bval_{(\orbit, C)}$
as promised above.

\begin{proposition} \label{indep}
The duality map is well-defined; 
that is, it is independent of the choices made for $s$ and $e$.
\end{proposition}

\begin{proof}

We show that the representation $E$ constructed above
is independent of $s$ (it is clearly
independent of $e$, since Springer representations only depend on 
the orbit through $e$).  

One of the main properties of Lusztig-Spaltenstein duality is that 
if $e \in \levi'$ where $\levi'$ is a Levi subalgebra of $\g$, 
then $$d_{\g}(\orbit^{\g}_e) = 
\mbox{Ind}_{\levi'}^{\g} d_{\levi'}(\orbit^{\levi'}_e).$$  
The notation on the right-side
is Lusztig-Spaltenstein induction.  
According to \cite{lusp} and the validity
of property (B), one therefore has that $$E_{d_{\g}(\orbit_e),1} = 
j_{W(\levi')}^W (E_{d_{\levi}(\orbit^{\levi}_e), 1})$$
where $W(\levi')$ is the Weyl group of $\levi'$.

Now let $S$ be a maximal torus in $Z_G(s, e)$.
Then $\levi' = Z_{\g}(s, S)$ is a Levi subalgebra of $\levi = Z_{\g}(s)$
and $e \in \levi' \subset \levi$. 
Moreover $\levi'$ is of the form $Z_{\g}(s')$ for a semisimple
element $s' \in Z_G(e)$ and the image of $s'$ in $A(e)$ 
necessarily belongs to $C$. 
By the transitivity of truncated 
induction applied to the
sequence $\levi' \subset \levi \subset \g$,
we see that the representation $E$
is the same whether we work with respect to $\levi$ or $\levi'$, i.e.
whether we work with $s$ or $s'$.
The main result of \cite{so1} is that the pair $(\levi', e)$ (up to simultaneous 
conjugation by elements in $G$) is determined by 
(and determines in the case when $G$ is of adjoint type) 
the conjugacy class $C \subset A(\orbit_e)$.  In other words, 
$(\levi', e)$ (up to simultaneous 
conjugation by elements in $G$)
is independent of the choice of $s$ and 
thus the construction of $E$ is independent of the choices made for
$s$ and $e$.

\end{proof}

Proposition \ref{keyprop} now follows from the above proof either 
by invoking property (B) in $L$, 
or more simply by applying the dimension formula
$\dim d_{\levi}(\orbit^{\levi}_e) = \dim(\levi) - \dim(\levi') + 
\dim d_{\levi'}(\orbit^{\levi'}_e)$ 
for induction of orbits from the Levi subalgebra $\levi'$ of $\levi$.

\section{Duality in Classical Groups}

We now calculate the duality map for the classical groups not of type $A$
(there is nothing new here in type $A$).
First, we assign to a partition 
$\lambda$ (which may or may not correspond to 
a nilpotent orbit in that classical group) a representation $E_{\lambda}$ 
of the Weyl group of that group.   
When $\lambda$ corresponds to 
an actual nilpotent orbit, $E_{\lambda}$ is just the Springer representation 
$E_{\lambda, 1}$
for the orbit with trivial local system.

Let $\lambda \in \parti(m)$ 
where $m$ is even in types $C,D$ and odd in type $B$.
Form the dual partition $\lambda^*$.  
Separate the parts of $\lambda^*$
into its odd parts $2 \alpha_1 + 1 \geq 2 \alpha_2 + 1 \geq 
2\alpha_3 + 1 \dots 
\geq 2 \alpha_r + 1$  
and its even parts $2 \beta_1 \geq 2 \beta_2 \geq 2 \beta_3 \dots
\geq 2 \beta_s$.
We then associate to $\lambda$ the representation 
$E_{\lambda} = j_{W'}^{W}(\epsilon_{W'})$
where $W'$ is the Weyl group of the subsystem of type
\begin{equation} \label{springy}
\Phi'= \begin{cases}
B_{\alpha_1} \times B_{\beta_1} \times D_{\alpha_2 + 1} \times D_{\beta_2}
\times B_{\alpha_3} \times B_{\beta_3} \times D_{\alpha_4 + 1} \times D_{\beta_4} \dots
& \text{in type $B$,} \\
D_{\alpha_1 + 1} \times C_{\beta_1} \times C_{\alpha_2} \times D_{\beta_2}
\times D_{\alpha_3+1} \times C_{\beta_3} \times C_{\alpha_4} \times D_{\beta_4} \dots
& \text{in type $C$,} \\
D_{\alpha_1 + 1} \times D_{\beta_1} \times B_{\alpha_2} \times B_{\beta_2}
\times D_{\alpha_3+1} \times D_{\beta_3} \times B_{\alpha_4} \times B_{\beta_4} \dots
& \text{in type $D$.} 
\end{cases}
\end{equation}
and $\epsilon_{W'}$ is the sign representation of $W' = W(\Phi')$.
In type $D$ we are thinking 
of the representation as truncated induction in type $B$
followed by restriction to $W(D_n) \subset W(B_n)$ which is known
to produce an irreducible representation as long as
$\lambda$ is not very even (we ignore the case where $\lambda$ is very even
since we never need to consider it in what follows).  

\begin{lemma} \label{springer}
When $\lambda$ corresponds to a nilpotent orbit in the appropriate 
classical group, 
$E_{\lambda}$ is the Springer representation of 
this orbit with the trivial local system.
\end{lemma}

\begin{proof}
This is shown by following Lusztig's version
of Shoji's algorithm (see chapter 13.3 of \cite{ca}).
\end{proof}

\begin{lemma} \label{collapse}
For two partitions $\lambda$, $\mu$  
with the same $X$-collapse where $X = B,C,$ or $D$,
we have $E_{\lambda} = E_{\mu}$ as representations of $W(X_n)$.
\end{lemma}

\begin{proof}
We give the proof in type $B$, the other types being similar.
Assume $\lambda \in \parti(2n+1)$, but $\lambda \notin \parti_B(2n+1)$.
List all the even parts of $\lambda$ in decreasing order as
$\lambda_{e_1} \geq \lambda_{e_2} \geq \dots \geq \lambda_{e_l}$.
Then there exists an $m$ such that $\lambda_{e_{2m-1}} > \lambda_{e_{2m}}$
since $\lambda \notin \parti_B(2n+1)$.  
Let $\mu$ be the partition obtained from $\lambda$ by 
replacing the part $\lambda_{e_{2m-1}}$ by $\lambda_{e_{2m-1}}-1$ 
and the part $\lambda_{e_{2m}}$ by $\lambda_{e_{2m}} + 1$
and leaving all other parts of $\lambda$ unchanged.
This is the basic $B$-collapsing move (see \cite{cm})
and it suffices to show that $E_{\mu} = E_{\lambda}$.
The dual partition $\mu^*$ equals $\lambda^*$ except
that $\mu^*_{\lambda_{e_{2m}}+1} = \lambda^*_{\lambda_{e_{2m}}+1}+1$
and $\mu^*_{\lambda_{e_{2m-1}}}  = \lambda^*_{\lambda_{e_{2m-1}}}-1$.  
Now write $\lambda^* = \alpha \cup \beta$ where $\alpha$ and $\beta$
consist respectively of the odd and even parts of $\lambda^*$.
Given a part $\lambda^*_i$ of $\lambda^*$, 
recall that we are calling
$i$ the position of $\lambda^*_i$ in $\lambda^*$.
Now $\lambda^*_i$ will occur as a part of either $\alpha$ or $\beta$.
We refer to its position in whichever partition $\alpha$ or $\beta$
it occurs in as its parity position.  We use similar language for $\mu^*$.

Assume first  
that $a = \lambda^*_{\lambda_{e_{2m}}+1}$ is even.
Then it is possible to show that the parity
position of $\lambda^*_{\lambda_{e_{2m}}+1}$ is odd
and the parity position of 
$\mu^*_{\lambda_{e_{2m}}+1} =
\lambda^*_{\lambda_{e_{2m}}+1}+1$ is also odd.
These parts will each 
contribute a root subsystem of type $B_{a/2}$ in  
the definition of $E_{\lambda}$ or $E_{\mu}$, respectively.
On the other hand, if $a$ is odd the parity position of these parts
is both even and they will each contribute 
a root subsystem of type $D_{(a+1)/2}$ in the definition of $E_{\lambda}$ 
or $E_{\mu}$, respectively.

Since there are no even parts of $\lambda$ between $\lambda_{e_{2m-1}}$
and $\lambda_{e_{2m}}$, we have 
$\lambda^*_{\lambda_{e_{2m}}+2j} = \lambda^*_{\lambda_{e_{2m}}+2j+1}$ 
for $j= 1, 2 , \dots , \frac{\lambda_{e_{2m-1}} - \lambda_{e_{2m}}}{2} - 1$.
Since the same equalities hold for $\mu^*$, we see that these parts
(which come in pairs) contribute the same terms to $E_{\lambda}$ and $E_{\mu}$.

Finally consider $b = \lambda^*_{\lambda_{e_{2m-1}}}$.  
If $b$ is even, then the parity position 
of both $\lambda^*_{\lambda_{e_{2m-1}}}$ and 
$\mu^*_{\lambda_{e_{2m-1}}}$ is even, so they contribute
a term $D_{b/2}$ in the definition of $E_{\lambda}$ 
or $E_{\mu}$, respectively.
And if $b$ is odd, then they both have odd parity position 
and their contribution is $B_{(b-1)/2}$.

Hence a basic collapsing move does not affect
the attached representation and the result is proved.
\end{proof}

\begin{remark}
The formulas for the $\bval$-values in the classical
groups are a consequence of the previous
propositions and the validity of property (B). 
\end{remark}

We now explain the bijection between $\parti^{sp}_B (2n+1)$
and $\parti^{sp}_C (2n)$ explicitly in terms of partitions
(see \cite{spa1}, \cite{kp}).  
Given $\lambda =[\lambda_1 \geq \dots 
\geq \lambda_{k-1} \geq \lambda_k] \in \parti_B (2n+1)$,
let $\lambda^{-} = [\lambda_1 \geq \dots \geq \lambda_{k-1}
\geq \lambda_k -1 ]$ 
and set $\lambda^C = \lambda^{-} _C$.
Then $\lambda^C \in \parti^{sp}_C(2n)$.
Similarly given $\lambda =[\lambda_1 \geq \dots 
\geq \lambda_k] \in \parti_C (2n)$, let 
$\lambda^{+} = [\lambda_1 +1 \geq \dots \geq \lambda_{k-1}
\geq \lambda_k]$
and
$\lambda_{+} = [\lambda_1 \geq \dots \geq \lambda_{k-1}
\geq \lambda_k \geq 1 ]$
and set $\lambda^B = (\lambda^{+})_{B}$.
Then $\lambda^B \in \parti^{sp}_B (2n+1)$ 
and moreover, $\lambda^B = (\lambda_{+})^{*}_B)^{*}$.
For $\lambda \in \parti_B (2n+1)$, we
have $(\lambda^C)^B = (\lambda^{*}_B)^{*}$; in particular,
if $\lambda$ is special, $(\lambda^C)^B = \lambda$.
Similarly for $\lambda \in \parti_C (2n)$, we
have $(\lambda^B)^C = (\lambda^{*}_C)^{*}$; in particular,
if $\lambda$ is special, $(\lambda^B)^C = \lambda$.

In what follows we identify representations of 
$W(B_n)$ and $W(C_n)$ via the isomorphism of these
two Coxeter groups which corresponds to
interchanging long and short roots.  

\begin{lemma} \label{transfer}
For $\lambda \in \parti_C (2n)$ we have 
$E_{\lambda^*} = E_{(\lambda^B)^*}$,
where the left side of the identity is computed in type $C$
and the right in type $B$.
For $\lambda \in \parti_B (2n+1)$ 
we have $E_{\lambda^*} = E_{(\lambda^C)^*}$,
where the left side of the identity is computed in type $B$
and the right in type $C$.
\end{lemma}

\begin{proof}
We prove the first isomorphism.
We noted above that $(\lambda^B)^* = (\lambda_{+})^*_B$
and so in type $B$,
$E_{(\lambda^B)^*} = E_{(\lambda_{+})^*}$ since by
the previous proposition we can omit the $B$-collapse on the right side.
To prove the desired identity
we must study the odd and even parts of 
$\lambda$ (in type $C$) and 
$\lambda_{+}$ (in type $B$).  
These partitions are the same except that the latter
has an extra part equal to $1$ at the end.  
Now because $\lambda \in \parti_C (2n)$ the definition
of the subsystem in equation (\ref{springy})
for $E_{\lambda^*}$ in type $C$ and 
for $E_{(\lambda_{+})^*}$ in type $B$ coincide
(with the extra part in $\lambda_{+}$ playing no role at all).

We now prove the second isomorphism.
The first isomorphism implies that
$E_{(\lambda^C)^*} = E_{((\lambda^C)^B)^*}$,
where the left side is in type $C$ and the right in type $B$.
Since $\lambda \in \parti_B (2n+1)$, we have
$(\lambda^C)^B = (\lambda^{*}_B)^{*}$ and thus 
$E_{((\lambda^C)^B)^*} = E_{(\lambda^{*}_B)} = E_{\lambda^{*}}$.
The last equality holds since in type $B$ we can omit the $B$-collapse.
\end{proof}

\begin{theorem} \label{recipe}
Our duality map $d: \nilcone_{o,c} \to {^L \nilcone_{o}}$ sends the pair 
$(\nu, \eta)$ to the orbit $\lambda$ according to 
the following recipe:

\begin{equation} \lambda = \begin{cases}
(\nu \cup \eta^C)^{*}_C & \text{$\g$ type $B$} \\
(\nu \cup \eta^B)^*_B & \text{$\g$ type $C$} \\
(\nu \cup (\eta^*_D)^*)^{*}_D & \text{$\g$ type $D$} \\
\end{cases}
\end{equation}
\end{theorem}

Note that the case of $\nu$ equal to the empty partition corresponds
to Lusztig-Spaltenstein duality.  In type $D$ if $\nu$ is non-empty,
our assumptions about 
$(\nu, \eta)$ ensure that $\eta$ is not very even. 

\begin{proof}
We may choose $s \in Z_G(e)$ representing $C$ so that 
$\levi = Z_{\g}(s)$ has semisimple rank equal to the rank of $\g$.
Then $e \in \levi$ is specified by
the pair of partitions $(\nu, \eta)$.  
Our first step is to compute the Springer representation $E_{d_{\levi}(e), 1}$ 
of $W(s)$ associated to $d_{\levi}(e)$.
The pair of partitions associated to $d_{\levi}(e)$ is 
\begin{equation} 
\begin{aligned}
( \nu^*_D, \eta^*_B ) & \text{ \ $\g$ type $B$} \\
( \nu^*_C, \eta^*_C ) & \text{ \ $\g$ type $C$} \\
( \nu^*_D, \eta^*_D ) & \text{ \ $\g$ type $D$} \\
\end{aligned}
\end{equation}
By lemmas \ref{springer} and \ref{collapse}, the associated Springer
representation $E_{d_{\levi}(e), 1}$ of  $W(s)$ 
is $E_{\nu^*} \boxtimes E_{\eta^*}$.
Consider this now as a representation of 
$W(D_k) \times W(C_{n-k})$ in type $B$, 
$W(C_k) \times W(B_{n-k})$ in type $C$, 
and $W(D_k) \times W(D_{n-k})$ in type $D$
(there is no change in type $D$).
Then by applying lemma \ref{transfer}
in types $B$ and $C$ and lemma \ref{collapse} again in type $D$, 
this representation can be described as
\begin{equation} \begin{aligned}
E_{\nu^*} \boxtimes E_{(\eta^C)^*} & \text{ \ $\g$ type $B$} \\
E_{\nu^*} \boxtimes E_{(\eta^B)^*} & \text{ \ $\g$ type $C$} \\
E_{\nu^*} \boxtimes E_{(\eta^*)_D} & \text{ \ $\g$ type $D$} 
\end{aligned}
\end{equation}
These representations possess property (B) as they possess property (B)
in each simple component.  Hence we
can apply truncated induction up to $W(\dualg)$.
Then by transitivity of induction we claim that we arrive at the representation
\begin{equation} \begin{aligned} \label{finish-up}
E_{(\nu \cup \eta^C)^*} & \text{ \ $\g$ type $B$} \\
E_{(\nu \cup \eta^B)^*} & \text{ \ $\g$ type $C$} \\
E_{(\nu \cup (\eta^*_D)^*)^*} & \text{ \ $\g$ type $D$} 
\end{aligned}
\end{equation}
where the first is a representation of $W(C_n)$, the second of $W(B_n)$,
and the third of $W(D_n)$.
In type $D$, we use the fact that $(\eta^*_D)^*$ belongs to $\parti_C(2n-k)$.
Therefore in all types if the multiplicity of $i$ in $\nu$ is odd,
then the multiplicity of $i$ in $\eta^C, \eta^B,$ or   
$(\eta^*_D)^*$, respectively, is even.
Then the validity of equation (\ref{finish-up}) 
follows from the definition of $E$ in equation (\ref{springy}). 

The proof is completed by applying lemma \ref{collapse} in $\dualg$.
\end{proof}

\section{Bookends}

\begin{proposition} \label{surjective}
Our duality map is surjective.
\end{proposition}

\begin{proof}
We verified this case-by-case (Lusztig already did this in
his work with his original map although the details are not recorded).
In fact, we will try to exhibit canonical elements $(\orbit, C)$ 
of $\nilcone_{o,c}$ which map bijectively to $^L \nilcone^{sp}_o$.
These are denoted by a star (*) in the tables for the exceptional groups
and we now explain their construction in the classical groups.

Assume $\dualg$ is of type $X$ where $X = B,C,D$
and $\lambda \in \parti_X(m)$ where $m$ is even or odd depending on $X$.
Consider $\lambda^*$.  We ask whether $\lambda^*$ belongs
to $\parti_B(m), \parti_C(m)$, or $\parti_C(m)$ depending on 
whether $X$ is $B$, $C$, or $D$, respectively.
In other words, we ask whether $\lambda$ is special (note the funny
situation in type $D$).
If not, we may uniquely write $\lambda^* = \nu \cup \mu$ where
$\nu$ is distinguished of type $C$, $D$, or $D$,
and where $\mu$ belongs to $\parti_B(m'), \parti_C(m')$, or $\parti_C(m')$,
for some $m'$, depending on whether $X$ is $B$, $C$, or $D$.

We now show that $\mu$ has the property that 
$\mu^* \in \parti_X(m')$.
This is because $\lambda \in \parti_X(m)$
and the process of forming $\mu^*$ amounts to 
taking $\lambda$ and diminishing some of its parts;
however,
parts not congruent to $\epsilon$ will only be diminished by an even number,
so the resulting partition remains of the same type
as $\lambda$.   Hence $\mu^* \in \parti_X(m')$ for some $m'$.
It follows that $\mu^* \in \parti^{sp}_X(m')$ since $\mu$ belonged to 
$\parti_B(m'), \parti_C(m')$, or $\parti_C(m')$ depending on $X$.
It is also true that $\mu$ is itself special in types $B$ and $C$.

We can now define $\eta$.
In type $B$ we set $\eta = \mu^C$;
in type $C$ we set $\eta = \mu^B$;
and in type $D$ we set $\eta = \mu_D$.
Now because $\mu^*$ is special (and $\mu$ is special in types 
$B$ and $C$), we have
$\eta ^B = (\mu^C)^B = \mu$ in type $B$;
$\eta ^C = (\mu^B)^C = \mu$  in type $C$; 
and  $(\eta^*_D)^* = (((\mu^*)^*_D)^*_D)^* = (\mu^*)^* = \mu$ in type $D$.
The second equality in type $D$ follows since
applying Lusztig-Spaltenstein duality twice (to $\mu^*$ in this case) 
is the identity on special orbits.

Thus in all types 
$d_{(\nu, \eta)} = (\nu \cup \mu)^*_X = (\lambda^*)^*_X 
= \lambda_X = \lambda $ 
where the last equality holds since $\lambda \in \parti_X(m)$. 
We conclude in all types that $(\nu, \eta)$ has the property
that $d_{(\nu, \eta)} = \lambda$.
\end{proof}

\begin{remark}
These canonical elements $(\orbit, C)$ that we have listed
which surject onto $^L \nilcone_o$  
have the property that $\orbit$ is always special.
In fact, all orbits $\orbit'$ of $^L \nilcone_o$ 
in the same special piece
are affiliated with the same special orbit $\orbit$ of $\g$
and in fact $\orbit = d_{(\orbit', 1)}$ 
(nilpotent orbits
are in the same special piece exactly when their dual orbits
are the same).
Hence for each orbit $\orbit'$ in $^L \nilcone_o$ we
get a conjugacy class in $\bar{A}(\orbit)$ where $\orbit = d_{(\orbit', 1)}$.
This should be the same conjugacy class that Lusztig attaches
to orbits in \cite{lu2}
under the identification  
$\bar{A}(\orbit) \cong \bar{A}(d_{(\orbit, 1)})$.
\end{remark}

\begin{proposition}
Let $C, C'$ be conjugacy classes in $A(\orbit)$ whose image
in $\bar{A}(\orbit)$ coincide, then 
$d_{(\orbit , C)}= d_{(\orbit , C')}$.
\end{proposition}

\begin{proof}
Again we checked this on a case-by-case basis. In the exceptional
groups, this amounts to a quick glance at the tables which follow.
In the classical groups, 
it requires attention to the computations
in the proof of theorem 
\ref{recipe}, together with the explicit description of the
canonical quotient.  
We omit the details.
\end{proof}


\newpage

\section{Duality in Exceptional Groups}

To compute the duality in the exceptional groups we used
knowledge of Lusztig-Spaltenstein duality and the Springer correspondence
(see \cite{ca}).  Furthermore, we 
computed truncated induction by using the induce/restrict tables of Alvis \cite{al}.
We have listed only those orbits with non-trivial component groups
when $G$ is of adjoint type since these are the only
orbits for which we are saying something new.  The stars(*) refer to
the (putatively) canonical pair $(\orbit, C)$ which maps to a given
orbit in
the dual Lie algebra as is done in proposition \ref{surjective} for the
classical groups.  

\begin{center}
\begin{tabular}{|c|c|c|c|c|} \hline
\multicolumn{5}{|c|}{$G_{2}$} \\ \hline
\multicolumn{1}{|r|}{$\DynkinG{.2}{}{}{}{}$} 
& \multicolumn{1}{|c|}{$(\levi, e)$}
& \multicolumn{1}{|c|}{$\bval$}
& \multicolumn{1}{|c|}{$\bar{A}(\orbit)$}
& \multicolumn{1}{|c|}{Dual} 

\\ \hline

\weightDynkinG{.2}{2}{0} & * $G_{2}(a_1)$ & 1 & $S_3$ & $G_{2}(a_1)$ \\ 
 & * $A_1 + \Tilde{A_{1}}$ & 2 & & $\Tilde{A_{1}}$ \\ 
 & * $A_2$ & 3 & & $A_{1}$ \\ \hline

\end{tabular}
\end{center}

\bigskip

\begin{center}
\begin{tabular}{|c|c|c|c|c|} \hline
\multicolumn{5}{|c|}{$F_{4}$} \\ \hline
\multicolumn{1}{|r|}{$\DynkinF{.2}{}{}{}{}$} 
& \multicolumn{1}{|c|}{$(\levi, e)$}
& \multicolumn{1}{|c|}{$\bval$}
& \multicolumn{1}{|c|}{$\bar{A}(\orbit)$}
& \multicolumn{1}{|c|}{Dual} 

\\ \hline


\weightDynkinF{.2}{0}{0}{0}{1} & * $\Tilde A_{1}$ & $1$ & $S_2$ & $F_{4}(a_{1})$ \\ 
    & $2A_{1}$ & 2 && $F_{4}(a_{2})$   \\ \hline


\weightDynkinF{.2}{2}{0}{0}{0} & * $A_{2}$ & 3 & $1$ & $B_3$ \\ 
    & $2A_{1} + \Tilde{A_{1}}$ & 3 & & $B_3$ \\ \hline



\weightDynkinF{.2}{2}{0}{0}{1} &  $B_{2}$ & 4 & $S_2$ & $F_{4}(a_{3})$ \\ 
  & $A_3$ & 6 & & $B_2$ \\ \hline


\weightDynkinF{.2}{1}{0}{1}{0} & $C_{3}(a_{1})$ & 4 &  $S_2$ & $F_{4}(a_{3})$  \\ 
   & $A_{1} + B_{2}$ & 5 & & $C_3(a_1)$ \\ \hline

\weightDynkinF{.2}{0}{2}{0}{0} & * $F_{4}(a_{3})$ & 4 &  $S_4$ & $F_{4}(a_{3})$ \\ 
        & * $A_{1} + C_{3}(a_{1})$ & 5 & & $C_{3}(a_{1})$\\  
        & * $A_{2} + \Tilde{A_{2}}$ & 6 & &$A_{1} + \Tilde{A_{2}}$\\ 
        & * $B_{4}(a_{2})$ & 6 & & $B_{2}$ \\ 
        & * $A_{3} + \Tilde{A_{1}}$ &  7 &  & $A_{2} + \Tilde{A_{1}}$ \\ \hline



\weightDynkinF{.2}{0}{2}{0}{2} & * $F_{4}(a_{2})$ & 10 & $1$ & $A_{1} + \Tilde{A_{1}}$ \\ 
                & $A_{1} + C_{3}$ & 10 & &  $A_{1} + \Tilde{A_{1}}$ \\ \hline 

\weightDynkinF{.2}{2}{2}{0}{2} & * $F_{4}(a_{1})$ & 13 & $S_2$ & $\Tilde{A_{1}}$ \\ 
                & * $B_{4}$ & 16 &  & $A_1$ \\ \hline


\end{tabular}
\end{center}
\bigskip
\begin{center}
\begin{tabular}{|c|c|c|c|c|} \hline
\multicolumn{5}{|c|}{$E_{6}$} \\ \hline
\multicolumn{1}{|c|}{\DynkinEsix{.18}} 
& \multicolumn{1}{|c|}{$(\levi, e)$}
& \multicolumn{1}{|c|}{$\bval$}
& \multicolumn{1}{|c|}{$\bar{A}(\orbit)$}
& \multicolumn{1}{|c|}{Dual} 
\\ \hline

\weightDynkinEsix{.18}{0}{0}{0}{0}{0}{2} & * $A_{2}$ & 3 & $S_2$ & $E_{6}(a_{3})$ \\ 
                 & * $4 A_{1}$ & 4 & & $A_{5}$\\ \hline

\weightDynkinEsix{.18}{0}{0}{2}{0}{0}{0} & * $D_{4}(a_{1})$  & 7 &  $S_3$ & $D_{4}(a_{1})$  \\ 
                 & * $A_{3} + 2A_{1}$ & 8 & &  $A_3 + A_1$ \\
                 & * $3A_{2}$ & 9 & & $2A_2 + A_1$ \\ \hline

\weightDynkinEsix{.18}{2}{0}{2}{0}{2}{0} & * $E_{6}(a_{3}) $ & 15  & $S_2$ & $A_{2}$ \\ 
                 & * $A_{5} + A_{1}$ & 16 & & $3A_1$ \\ \hline
\end{tabular}
\end{center}
\bigskip

\newpage
\begin{center}
\begin{tabular}{|c|c|c|c|c|} \hline
\multicolumn{5}{|c|}{$E_{7}$} \\ \hline
\multicolumn{1}{|c|}{\DynkinEseven{.18}} 
& \multicolumn{1}{|c|}{$(\levi, e)$}
& \multicolumn{1}{|c|}{$\bval$} 
& \multicolumn{1}{|c|}{$\bar{A}(\orbit)$}
& \multicolumn{1}{|c|}{Dual} 
\\ \hline
\weightDynkinEseven{.18}{2}{0}{0}{0}{0}{0}{0} & * $A_{2}$ & 3 & $S_2$ & $E_7(a_3)$ \\ 
                 & * $(4 A_{1})'$ & 4 & & $D_6$ \\ \hline

\weightDynkinEseven{.18}{1}{0}{0}{0}{1}{0}{0} & * $A_{2}+A_{1}$ & 4 & $S_2$ & $E_6(a_1)$\\ 
                 & $5 A_{1}$ & 5 & & $E_7(a_4)$  \\ \hline

\weightDynkinEseven{.18}{0}{2}{0}{0}{0}{0}{0} & * $D_{4}(a_{1})$  & 7 & $S_3$ & $E_7(a_5)$\\ 
                 & * $3A_{2}$ & 9 & & $A_5+A_1$ \\ 
                 & * $(A_{3} + 2A_{1})'$ & 8 & & $D_6(a_2)$ \\ \hline

\weightDynkinEseven{.18}{0}{1}{0}{0}{0}{1}{1} & * $D_{4}(a_{1})+A_{1} $ & 8 & $S_2$ & $E_6(a_3)$ \\ 
                 & * $A_{3} + 3A_{1}$ & 9 & & $(A_5)'$ \\ \hline

\weightDynkinEseven{.18}{0}{0}{1}{0}{1}{0}{0} & * $A_{3}+A_{2} $ & 9 & $1$ & $D_5(a_1)+A_1$ \\ 
                 & $D_{4}(a_{1}) + 2A_{1}$ & 9 & & $D_5(a_1)+A_1$ \\ \hline

\weightDynkinEseven{.18}{2}{0}{0}{0}{2}{0}{0} & * $A_{4}$ & 10 & $S_2$ & $D_5(a_1)$  \\ 
                 & * $2A_{3}$ & 12 & &  $D_4+A_1$\\ \hline

\weightDynkinEseven{.18}{1}{0}{1}{0}{1}{0}{0} & * $A_{4} + A_{1}$ & 11 & $S_2$ & $A_4+A_1$ \\ 
                 & $A_{1}+2A_{3}$ & 13 & & $A_3+A_2+A_1$ \\ \hline

\weightDynkinEseven{.18}{2}{0}{1}{0}{1}{0}{0} & * $D_{5}(a_{1})$ & 13 & $S_2$ & $A_4$\\ 
                 & $D_{4} + 2A_{1}$ & 14 & & $A_3+A_2$ \\ \hline

\weightDynkinEseven{.18}{0}{2}{0}{0}{2}{0}{0} & * $E_{6}(a_{3})$ & 15 &  $S_2$ & $D_4(a_1)+A_1$  \\
                 & * $(A_{5} + A_{1})'$ & 16 & & $A_3+2A_1$ \\ \hline

\weightDynkinEseven{.18}{0}{0}{2}{0}{0}{2}{0} & * $E_{7}(a_{5})$ & 16 & $S_3$ & $D_4(a_1)$ \\ 
                 & * $A_{5}+A_{2}$ & 18 &  & $2A_2+A_1$ \\ 
                 & * $A_{1}+D_{6}(a_{2})$ & 17 & & $(A_3+A_1)'$ \\ \hline

\weightDynkinEseven{.18}{2}{0}{2}{0}{0}{2}{0} & * $E_{7}(a_{4})$ & 22 & $1$ & $A_2+2A_1$ \\ 
                 & $A_{1}+D_{6}(a_{1})$ &  22 &  & $A_2+2A_1$ \\ \hline

\weightDynkinEseven{.18}{2}{0}{2}{0}{2}{0}{0} & * $E_{6}(a_{1})$ & 25 & $S_2$ & $A_2+A_1$\\ 
                 &  * $A_{7}$ & 28 & & $4A_1$ \\ \hline

\weightDynkinEseven{.18}{2}{0}{2}{0}{2}{2}{0} & * $E_{7}(a_{3})$ & 30 & $S_2$ & $A_2$ \\ 
                 &  * $A_{1}+D_{6}$ & 31 & &  $(3A_1)'$ \\ \hline

\end{tabular}
\end{center}
\begin{center}
\begin{tabular}{|c|c|c|c|c|} \hline
\multicolumn{5}{|c|}{$E_{8}$} \\ \hline
\multicolumn{1}{|c|}{\DynkinEeight{.18}} 
& \multicolumn{1}{|c|}{$(\levi, e)$}
& \multicolumn{1}{|c|}{$\bval$} 
& \multicolumn{1}{|c|}{$\bar{A}(\orbit)$}
& \multicolumn{1}{|c|}{Dual} 
\\ \hline
\weightDynkinEeight{.18}{0}{0}{0}{0}{0}{0}{2}{0} & * $A_{2}$ & 3 & $S_2$ & $E_8(a_3)$\\ 
                 & * $(4 A_{1})''$ & 4 & & $E_7$ \\ \hline

\weightDynkinEeight{.18}{1}{0}{0}{0}{0}{0}{1}{0} & * $A_{2}+A_{1}$ & 4  & $S_2$& $E_8(a_4)$ \\ 
                 & $5 A_{1}$ & 5 & & $E_8(b_4)$ \\ \hline

\weightDynkinEeight{.18}{2}{0}{0}{0}{0}{0}{0}{0} & * $2A_{2}$ & 6  & $S_2$& $E_8(a_5)$ \\ 
                 & * $A_{2}+4 A_{1}$ & 7 & & $D_7$ \\ \hline

\weightDynkinEeight{.18}{0}{0}{0}{0}{0}{2}{0}{0} & * $D_{4}(a_{1})$  & 7  & $S_3$& $E_8(b_5)$ \\ 
                 &  * $3A_{2}$ & 9 & & $E_6+A_1$\\ 
                 &  * $(A_{3} + 2A_{1})''$ & 8 & & $E_7(a_2)$ \\ \hline

\weightDynkinEeight{.18}{0}{0}{0}{0}{0}{1}{0}{1} & *  $D_{4}(a_{1})+A_{1} $ & 8 & $S_3$ & $E_8(a_6)$\\ 
                 &  $3A_{2} + A_{1}$ & 10 & & $E_8(b_6)$ \\ 
                 &  $A_{3} + 3A_{1}$ & 9 & & $D_7(a_1)$ \\ \hline

\weightDynkinEeight{.18}{1}{0}{0}{0}{1}{0}{0}{0} & * $A_{3}+A_{2}$ & 9 & $1$ & $D_7(a_1)$ \\ 
                 &  $D_{4}(a_{1}) + 2A_{1}$ & 9 &  & $D_7(a_1)$ \\ \hline

\weightDynkinEeight{.18}{0}{0}{0}{0}{0}{0}{0}{2} & * $D_{4}(a_{1}) + A_{2}$ & 10 & $S_2$ & $E_8(b_6)$ \\ 
                 &  * $A_{3}+A_{2}+2A_{1}$ & 11 & & $A_7$ \\ \hline

\weightDynkinEeight{.18}{2}{0}{0}{0}{0}{0}{2}{0} & *  $A_{4}$ & 10 & $S_2$ & $E_7(a_3)$ \\ 
                 &  * $(2A_{3})''$ & 12 & & $D_6$  \\ \hline

\weightDynkinEeight{.18}{1}{0}{0}{0}{1}{0}{1}{0} &  * $A_{4} + A_{1}$ & 11 & $S_2$ & $E_6(a_1)+A_1$ \\ 
                 &   $A_{1}+2A_{3}$ & 13 & & $D_5+A_2$ \\ \hline

\weightDynkinEeight{.18}{0}{0}{1}{0}{0}{0}{1}{0} & * $A_{4} + 2A_{1}$ & 12 & $S_2$& $D_7(a_2)$ \\ 
                 &   $D_{4}(a_{1}) + A_{3}$ & 13 & & $D_5+A_2$ \\ \hline

\weightDynkinEeight{.18}{1}{0}{0}{0}{1}{0}{2}{0} & * $D_{5}(a_{1})$ & 13  & $S_2$& $E_6(a_1)$ \\ 
                 &   $D_{4} + 2A_{1}$ & 14 & & $E_7(a_4)$ \\ \hline
 
\weightDynkinEeight{.18}{0}{0}{0}{0}{0}{0}{2}{2} & * $D_{4}+A_{2}$ & 15 & $1$ & $A_6$ \\ 
                 &   $D_{5}(a_{1}) + 2A_{1}$ & 15 & & $A_6$ \\ \hline

\weightDynkinEeight{.18}{2}{0}{0}{0}{0}{2}{0}{0} & * $E_{6}(a_{3})$ & 15 & $S_2$ & $D_6(a_1)$  \\ 
                &  *  $(A_{5} + A_{1})''$ & 16 & & $D_5+A_1$ \\ \hline

\weightDynkinEeight{.18}{0}{1}{0}{0}{0}{1}{0}{1} &  $D_{6}(a_{2})$ & 16 & $S_2$& $E_8(a_7)$ \\ 
                 &  $D_{4}+A_{3}$ & 18 & &  $D_{6}(a_{2})$\\ \hline

\weightDynkinEeight{.18}{1}{0}{0}{1}{0}{1}{0}{0} & $E_{6}(a_{3})+A_{1}$ & 16 & $S_2$ & $E_8(a_7)$ \\ 
                 &  $A_{5} + 2 A_{1}$ & 17 & &  $E_7(a_5)$ \\ \hline

\weightDynkinEeight{.18}{0}{0}{1}{0}{1}{0}{0}{0} & $E_{7}(a_{5})$  & 16 & $S_3$& $E_8(a_7)$ \\ 
                 &  $A_{5}+A_{2}$ & 18 & & $E_6(a_3)+A_1$ \\ 
                 &  $A_{1}+D_{6}(a_{2})$ & 17 & & $E_7(a_5)$ \\ \hline

\end{tabular}

\end{center}

\newpage
\begin{center}
\begin{tabular}{|c|c|c|c|c|} \hline
\multicolumn{5}{|c|}{$E_{8}$} \\ \hline
\multicolumn{1}{|c|}{\DynkinEeight{.18}} 
& \multicolumn{1}{|c|}{$(\levi, e)$}
& \multicolumn{1}{|c|}{$\bval$}
& \multicolumn{1}{|c|}{$\bar{A}(\orbit)$}
& \multicolumn{1}{|c|}{Dual} 
\\ \hline
\weightDynkinEeight{.18}{0}{0}{0}{2}{0}{0}{0}{0} & * $E_{8}(a_{7})$  & 16 & $S_5$ & $E_8(a_7)$ \\ 
                 &  * $A_{5}+A_{2}+A_{1}$ & 19 &  &$A_5+A_1$  \\ 
                 &  * $2A_{4}$ & 20 &  & $A_4+A_3$ \\ 
                 &  * $D_{5}(a_{1}) + A_{3}$ & 19  &  & $D_5(a_1)+A_2$  \\
                 &  * $D_{8}(a_{5})$ & 18  & & $D_6(a_2)$ \\
                 &  * $E_{7}(a_{5})+A_{1}$ & 17  && $E_7(a_5)$ \\
                 &  * $E_{6}(a_{3})+A_{2}$ & 18  &&   $E_6(a_3)+A_1$ \\ \hline

\weightDynkinEeight{.18}{0}{1}{0}{0}{0}{1}{2}{1} & * $D_{6}(a_{1})$ & 21 & $S_2$& $E_6(a_3)$ \\ 
                 &  * $D_{5}+2A_{1}$ & 22 & & $A_5$ \\ \hline

\weightDynkinEeight{.18}{0}{0}{1}{0}{1}{0}{2}{0} &  * $E_{7}(a_{4})$ & 22 & $1$ & $D_5(a_1)+A_1$ \\ 
                 &  $A_{1}+D_{6}(a_{1})$ &  22 &  &$D_5(a_1)+A_1$ \\ \hline

\weightDynkinEeight{.18}{0}{0}{0}{2}{0}{0}{2}{0} & * $D_{5}+A_{2}$ & 23 &  $1$ &$A_4+A_2$ \\
                 &  $E_{7}(a_{4})+A_{1}$ & 23 &   &$A_4+A_2$ \\ \hline

\weightDynkinEeight{.18}{1}{0}{1}{0}{1}{0}{1}{0}  & * $D_{7}(a_{2})$ & 24 &  $S_2$& $A_4+2A_1$ \\ 
                 &  * $D_{5}+A_{3}$ & 26 & & $2A_3$ \\ \hline

\weightDynkinEeight{.18}{2}{0}{0}{0}{2}{0}{2}{0} &  * $E_{6}(a_{1})$ & 25 & $S_2$& $D_5(a_1)$ \\ 
                 &  * $(A_{7})''$ & 28 &  &$D_4+A_1$ \\ \hline

\weightDynkinEeight{.18}{1}{0}{1}{0}{1}{0}{2}{0} & * $E_{6}(a_{1})+A_{1}$ & 26 & $S_2$& $A_4+A_1$ \\ 
                 &  $A_{7} + A_{1}$ & 29 &  &$A_3+A_2+A_1$ \\ \hline

\weightDynkinEeight{.18}{0}{0}{2}{0}{0}{0}{2}{0} & * $E_{8}(b_{6})$  & 28 & $S_2$& $D_4(a_1) +A_2$ \\ 
                 &  $E_{6}(a_{1})+A_{2}$ & 28 & & $D_4(a_1) +A_2$ \\ 
                 &  * $D_{8}(a_{3})$ & 29 & & $A_3+A_2+A_1$ \\ \hline

\weightDynkinEeight{.18}{2}{0}{1}{0}{1}{0}{2}{0} & * $E_{7}(a_{3})$ & 30 & $S_2$& $A_4$ \\ 
                 &  $A_{1}+D_{6}$ & 31 & & $A_3+A_2$ \\ \hline

\weightDynkinEeight{.18}{2}{0}{0}{2}{0}{0}{2}{0} & * $D_{7}(a_{1})$ & 31 &  $1$ & $A_3+A_2$ \\ 
                 &  $E_{7}(a_{3})+A_{1}$ & 31 & & $A_3+A_2$ \\ \hline

\weightDynkinEeight{.18}{0}{0}{2}{0}{0}{2}{0}{0} & * $E_{8}(a_{6})$  & 32 &  $S_3$ & $D_4(a_1)+A_1$ \\ 
                 &  * $A_{8}$ & 36 & &$2A_2+2A_1$ \\ 
                 &  * $D_{8}(a_{2})$ & 34 & & $A_3+2A_1$ \\ \hline

\weightDynkinEeight{.18}{0}{0}{2}{0}{0}{2}{2}{0} & * $E_{8}(b_{5})$  & 37 &  $S_3$&$D_4(a_1)$ \\ 
                 &  * $E_{6}+A_{2}$ & 39 & & $2A_2+A_1$\\ 
                 &  * $E_{7}(a_{2})+A_{1}$ & 38 &  & $A_3+A_1$ \\ \hline

\weightDynkinEeight{.18}{2}{0}{2}{0}{0}{2}{0}{0} & * $E_{8}(a_{5})$ & 42 & $S_2$& $2A_2$ \\ 
                 &  * $D_{8}(a_{1})$ & 43 &  & $A_2+3A_1$ \\ \hline

\weightDynkinEeight{.18}{2}{0}{2}{0}{0}{2}{2}{0} & * $E_{8}(b_{4})$ & 47 &  $1$ & $A_2+2A_1$ \\ 
                 &  $E_{7}(a_{1}) + A_{1}$ & 47 &  & $A_2+2A_1$ \\ \hline

\weightDynkinEeight{.18}{2}{0}{2}{0}{2}{0}{2}{0} & * $E_{8}(a_{4})$ & 52 &  $S_2$& $A_2+A_1$ \\ 
                 &  * $D_{8}$ & 56 &&  $4A_1$ \\ \hline

\weightDynkinEeight{.18}{2}{0}{2}{0}{2}{2}{2}{0} & * $E_{8}(a_{3})$ & 63 &  $S_2$& $A_2$  \\ 
                 &  * $E_{7}+A_{1}$ & 64 & & $3A_1$ \\ \hline

\end{tabular}

\end{center}


\begin{thebibliography}{df}

\bibitem[AL]{al} D. Alvis, \emph{Induce/restrict matrices for exceptional Weyl groups},
available at \texttt{www.iusb.edu/$\sim$dalvis/}

\bibitem[BV]{bv} D. Barbasch and D. Vogan, Jr., 
{\em  Unipotent representations of complex semisimple groups},
 Ann. of Math. (2), \textbf{121} (1985), no. 1, 41--110.
 
\bibitem[BM]{bm}
W. Borho and R. MacPherson, {\em Repr\'esentations des groupes de Weyl et homologie
d'intersection pour les vari\'etes nilpotentes},
 C. R. Acad. Sci. Paris S\'er. I Math., \textbf{292} (1981), no. 15, 707--710. 
\bibitem[Ca]{ca} R. W. Carter, \emph{Finite Groups of Lie Type.}
John Wiley and Sons: Chichester, 1985.

\bibitem[CM]{cm} D. Collingwood and W. McGovern, \emph{Nilpotent Orbits in Semisimple Lie algebras.} Van Nostrand Reinhold: New York, 1993.

\bibitem[KP]{kp} H. Kraft and C. Procesi, 
\emph{A special decomposition of the nilpotent cone
of a classical Lie algebra}, Ast\'erique \textbf{173-174} (1989), 271-279.
 
\bibitem[L1]{lu1} G. Lusztig, \emph{A class of irreducible representations of a Weyl group 
I, II},
Indag. Math. \textbf{41} (1979), no. 3, 323--335; \textbf{44} (1982), no. 2, 219--226.

\bibitem[L2]{lu3} G.~Lusztig, {\em Characters of reductive groups over a finite field}.
Annals of Mathematics Studies, \textbf{107}. Princeton University Press, 
Princeton, N.J., 1984.

\bibitem[L3]{lu2} G. Lusztig, \emph{Notes on unipotent classes},
Asian J. Math. \textbf{1} (1997), no. 1, 194--207.

\bibitem[LuSp]{lusp} G.~Lusztig and N. Spaltenstein, \emph{Induced unipotent classes},
J. London Math. Soc. (2), \textbf{19} (1979), no. 1, 41--52. 

\bibitem[So1]{so2} E. Sommers, 
\emph{Nilpotent orbits and the affine flag manifold}, 
Ph.D. Thesis, M.I.T., 1997.

\bibitem[So2]{so1} E. Sommers, \emph{A generalization of the Bala-Carter theorem}, IMRN 
(Intl. Math. Res. Notices), 1998, 539-562.

\bibitem[Sp]{spa1} N. Spaltenstein, \emph{Classes unipotentes et sous-groupes de Borel.}
Lecture Notes in Mathematics, \textbf{946}. Springer-Verlag, Berlin-New York, 1982.

\end{thebibliography}
\end{document}